%% file: 2020_07_heat_integration.tex
\documentclass[preprint]{elsarticle}

\input{scripts_els}

\journal{a journal}

\usepackage{todonotes}

\begin{document}

\begin{frontmatter}

\title{Semi-analytic integration for a parallel space-time \\ boundary element method modeling the heat equation}

\author[it4iaddress,dept470]{Jan Zapletal\texorpdfstring{\corref{mycorrespondingauthor}}{}}
\cortext[mycorrespondingauthor]{Corresponding author}
\ead{jan.zapletal@vsb.cz}

\author[tugraz]{Raphael Watschinger}

\author[tugraz]{Günther Of}

\author[it4iaddress,dept470]{Michal Merta}

\address[it4iaddress]{IT4Innovations, VŠB -- Technical University of Ostrava. \\ 17. listopadu 2172/15, 708 00 Ostrava-Poruba, Czech Republic}

\address[dept470]{Department of Applied Mathematics, VŠB -- Technical University of Ostrava. \\ 17. listopadu 2172/15, 708 00 Ostrava-Poruba, Czech Republic}

\address[tugraz]{Institute of Applied Mathematics, Graz University of Technology. \\ Steyrergasse 30, A-8010 Graz, Austria}

\begin{abstract}
  The presented paper concentrates on the boundary element method (BEM) for the heat equation in three spatial dimensions. In particular, we deal with tensor product space-time meshes allowing for quadrature schemes analytic in time and numerical in space. The spatial integrals can be treated by standard BEM techniques known from three dimensional stationary problems. The contribution of the paper is twofold. First, we provide temporal antiderivatives of the heat kernel necessary for the assembly of BEM matrices and the evaluation of the representation formula. Secondly, the presented approach has been implemented in a publicly available library besthea allowing researchers to reuse the formulae and BEM routines straightaway. The results are validated by numerical experiments in an HPC environment. 
\end{abstract}


\end{frontmatter}


\section{Introduction}
\label{sec:01_intro}
\noindent
\input{01_intro}

\section{Boundary integral equations}
\label{sec:02_bie}
\input{02_bie}

\section{Boundary element method}
\label{sec:03_bem}
\input{03_bem}

\section{Summary}
\label{sec:035_tldr}
\input{035_tldr}

\section{Implementation}
\label{sec:04_imp}
\input{04_imp}

\section{Numerical experiments}
\label{sec:05_exp}
\input{05_exp}

\section{Conclusion}
\label{sec:06_con}
\input{06_con}

\section*{Acknowledgements}
The authors acknowledge the support provided by the Czech Science Foundation under the project 19-29698L, by the Austrian Science Fund (FWF) under the project I 4033-N32, and by The Ministry of Education, Youth and Sports from the Large Infrastructures for Research, Experimental Development, and Innovations project e-INFRA CZ -- LM2018140.

\bibliography{references}

\end{document}

%% file: scripts_els.tex
\usepackage[utf8x]{inputenc}
\usepackage[T1]{fontenc}
\usepackage{amsmath,amssymb,amsthm}
\usepackage[left=34mm,right=34mm]{geometry}
\usepackage{fixmath}
\usepackage{bm}
\usepackage{xcolor}
\usepackage{booktabs}
\usepackage{listings}
\usepackage{lmodern}
\PassOptionsToPackage{hyphens}{url}\usepackage{hyperref}

\graphicspath{{fig/}}

\definecolor{codegreen}{rgb}{0,0.6,0}
\definecolor{codegray}{rgb}{0.5,0.5,0.5}
\definecolor{codepurple}{rgb}{0.58,0,0.82}
\definecolor{backcolour}{rgb}{1,1,1}
\definecolor{codeblue}{RGB}{0,160,238}
\definecolor{codeorange}{RGB}{217,83,25}
\lstdefinestyle{mystyle}{
    language=C++,
    backgroundcolor=\color{backcolour},   
    commentstyle=\color{codeblue},
    keywordstyle=\color{codeorange},
    numberstyle=\tiny\color{codegray},
    breakatwhitespace=false,         
    breaklines=true,                 
    captionpos=b,                    
    keepspaces=true,                 
    numbers=left,                    
    numbersep=5pt,                  
    showspaces=false,                
    showstringspaces=false,
    showtabs=false,                  
    tabsize=2,
    basicstyle={\scriptsize\ttfamily},
    xleftmargin=.05\textwidth, 
    xrightmargin=.05\textwidth,
    morekeywords={omp,offload,target,parallel,simd,critical,atomic,linear,declare,uniform,reduction,offload_transfer},
    framextopmargin=2pt,
    frame=tb,
    escapeinside={(|}{|)},
    moredelim=[is][\color{codeorange}]{|*}{*|},
    showlines=true
}
\newcommand{\vect}[1]{\ensuremath{{\bm{#1}}}}

\newcommand{\vg}{{\vect{g}}}
\newcommand{\vh}{{\vect{h}}}

\newcommand{\vn}{{\vect{n}}}

\newcommand{\vr}{{\vect{r}}}
\newcommand{\vs}{{\vect{s}}}

\newcommand{\vu}{{\vect{u}}}

\newcommand{\vw}{{\vect{w}}}
\newcommand{\vx}{{\vect{x}}}
\newcommand{\vy}{{\vect{y}}}
\newcommand{\vz}{{\vect{z}}}

\newcommand{\vF}{{\vect{F}}}

\newcommand{\vzero}{{\vect{0}}}

\newcommand{\mat}[1]{\ensuremath{{\mathsf{#1}}}}

\newcommand{\matA}{{\mat{A}}}

\newcommand{\matD}{{\mat{D}}}

\newcommand{\matK}{{\mat{K}}}

\newcommand{\matM}{{\mat{M}}}
\newcommand{\matO}{{\mat{O}}}

\newcommand{\matT}{{\mat{T}}}

\newcommand{\matV}{{\mat{V}}}

\newcommand{\dif}{\ensuremath{{\mathrm{d}}}}

\newcommand{\R}{{\mathbb{R}}}

\newcommand{\cpp}{C\nolinebreak\hspace{-.05em}\raisebox{.4ex}{\tiny\bf +}\nolinebreak\hspace{-.10em}\raisebox{.4ex}{\tiny\bf +}}

\DeclareMathOperator{\vspan}{span}

\DeclareMathOperator{\supp}{supp}
\DeclareMathOperator{\curl}{\mathbf{curl}}

\DeclareMathOperator{\erf}{erf}
\DeclareMathOperator{\argmin}{argmin}

\numberwithin{equation}{section}
\numberwithin{figure}{section}
\numberwithin{table}{section}
\AtBeginDocument{\counterwithin{lstlisting}{section}}

\theoremstyle{remark}

\newtheorem*{remark*}{Remark}

%% file: 01_intro.tex
For a bounded Lipschitz domain $\Omega \subset \R^3$ we aim to solve the heat equation
\begin{equation}
  \label{eq:heat_eq}
  \frac{\partial u}{\partial t}(\vx,t) - \alpha \upDelta u(\vx,t) = 0 \quad\text{for } (\vx,t) \in \Omega \times (0,T) =: Q
\end{equation}
with the heat capacity constant $\alpha > 0$, the initial condition
\begin{equation}\label{eq:initial:cond}
  u(\vx,0) = 0 \quad\text{for }\vx \in \Omega
\end{equation}
and a Dirichlet or Neumann boundary condition, i.e.
\begin{equation*}
  u(\vx,t) = g(\vx,t) \quad\text{for }(\vx,t) \in \partial\Omega \times (0,T) =: \Sigma
\end{equation*}
or
\begin{equation*}
  \alpha\frac{\partial u}{\partial \vn}(\vx,t) = h(\vx,t) \quad\text{for }(\vx,t) \in \Sigma,
\end{equation*}
respectively.

Such initial boundary value problems can be solved by boundary element methods. A survey on discretisation methods involving boundary integral equations is given in \cite{Cos2004}.
Here we consider a space-time formulation and a Galerkin method for discretisation. A comprehensive analysis of the involved integral equations is given in \cite{Cos1990}.  Error analysis for the Galerkin method has been
provided in \cite{Cos1990,Noo1988,ArnNoo1987,DohNiiSte2019,Doh2019}. 
A space-time formulation has certain advantages with respect to adaptivity and parallelisation. It allows quite general adaptivity in space and time compared to time stepping and convolution quadrature methods. A common parallelisation in space can be enhanced by an additional parallelisation with respect to time, which is not possible for time-stepping approaches.

We aim to provide a complete and (hopefully) error-free presentation of details on the implementation of a Galerkin boundary element method for the three-dimensional heat equation considering all boundary integral operators.
Galerkin methods have been considered in, e.g., \cite{Cos1990,Noo1988,Doh2019,DohEtAl2018} for 2d and \cite{Mes2014,MesSchTau2014,MesSchTau2015} for 3d. 
Typically, implementational aspects are discussed only briefly and a lot of effort is necessary to transform the theoretical results into a performant computer code. Our aim is to remove this setback by providing a detailed discussion and a publicly available \cpp{} library.  

In case of space-time tensor product discretisations the integrals with respect to time can be carried out analytically. This may result in a significant reduction of computational times. 
In \cite{Noo1988}, aspects of the temporal integration are discussed and \cite{Mes2014} considers the 3d setting. 
Unfortunately such presentations are typically very brief and we try to fill the gap by a detailed discussion.

The paper is organised as follows. In Section~\ref{sec:02_bie} we introduce the considered space-time boundary integral equations. The discretisation by the boundary element method is provided in Section~\ref{sec:03_bem} together with the derivation of heat kernel antiderivatives necessary for the assembly of system matrices and the evaluation of the representation formula. For a later reference by an interested reader using the presented results we provide a short summary of the formulae in Section~\ref{sec:035_tldr}. Section~\ref{sec:04_imp} describes the implementation approach as provided in the besthea \cpp{} library~\cite{besthea}. We validate the result by numerical experiments in Section~\ref{sec:05_exp} and conclude in Section~\ref{sec:06_con}.

%% file: 02_bie.tex
The solution to the initial problem~\eqref{eq:heat_eq}--\eqref{eq:initial:cond} is given by the representation formula
\begin{equation}
  \label{eq:repr}
  u(\vx,t) = \widetilde V \Big(\alpha\frac{\partial u}{\partial \vn}\Big)(\vx,t) - Wu(\vx,t)
\end{equation}
with the single-layer potential
\begin{equation*}
  \widetilde V \Big(\alpha\frac{\partial u}{\partial \vn}\Big) (\vx,t) := \int_0^t \int_{\partial\Omega} G_\alpha(\vx-\vy,t-\tau) \alpha\frac{\partial u}{\partial \vn} (\vy,\tau) \,\dif\vs_\vy \,\dif\tau,
\end{equation*}
the double-layer potential
\begin{equation*}
  W u (\vx,t) := \int_0^t \int_{\partial\Omega} \alpha\frac{\partial G_\alpha}{\partial \vn_\vy}(\vx-\vy,t-\tau) u(\vy,\tau) \,\dif\vs_\vy \,\dif\tau,
\end{equation*}
the fundamental solution to the heat equation
\begin{equation*}
  G_\alpha(\vx-\vy,t-\tau) :=
  \begin{cases}
  \displaystyle
  \frac{1}{(4\pi\alpha(t-\tau))^{3/2}}\exp\bigg( -\frac{\|\vx-\vy\|^2}{4\alpha(t-\tau)} \bigg) & \text{for } t > \tau, \\
  0 & \text{otherwise,}
  \end{cases}
\end{equation*}
and its scaled normal derivative
\begin{equation*}
  \alpha\frac{\partial G_\alpha}{\partial \vn_\vy}(\vx-\vy,t-\tau) :=
  \begin{cases}
  \displaystyle
  \frac{(\vx-\vy)\cdot \vn_\vy}{16(\pi\alpha)^{3/2}(t-\tau)^{5/2}}\exp\bigg( -\frac{\|\vx-\vy\|^2}{4\alpha(t-\tau)} \bigg) & \text{for } t > \tau, \\
  0 & \text{otherwise.}
  \end{cases}
\end{equation*}  
The operators $\widetilde{V}$ and $W$ are well-defined in the setting of anisotropic Sobolev spaces, see e.g.~\cite{LioMag1972Vol1,LioMag1972Vol2} for a definition of such spaces. In particular, natural choices are $X := H^{1/2,1/4}(\Sigma)$ for the space of the Dirichlet datum $u$ and its dual $X' := H^{-1/2,-1/4}(\Sigma)$ for the Neumann datum $w :=\alpha \partial u / \partial \vn$.

By applying the Dirichlet and Neumann trace operators to the representation formula \eqref{eq:repr} we get the boundary integral equations \cite{Cos2004,Cos1990,Noo1988,DohNiiSte2019}
\begin{align*}
  Vw(\vx,t) &= \bigg(\frac{1}{2} I + K \bigg) u(\vx,t) &&\text{for almost all } (\vx,t) \in \Sigma, \\
  Du(\vx,t) &= \bigg(\frac{1}{2} I - K_T' \bigg) w(\vx,t) &&\text{for almost all } (\vx,t) \in \Sigma,
\end{align*}
respectively. The boundary integral operators $V$, $K$, $D$, and $K_T'$ satisfy
\begin{align*}
  V &\colon X' \to X, &Vw(\vx,t) &= \int_0^t \int_{\partial\Omega} G_\alpha(\vx-\vy,t-\tau) w(\vy,\tau) \,\dif\vs_\vy \,\dif\tau, \\
  K &\colon X \to X, &Ku(\vx,t) &= \int_0^t \int_{\partial\Omega} \alpha\frac{\partial G_\alpha}{\partial \vn_\vy}(\vx-\vy,t-\tau) u(\vy,\tau) \,\dif\vs_\vy \,\dif\tau, \\
  D &\colon X \to X', &Du(\vx,t) &= -\alpha \frac{\partial}{\partial \vn_\vx} \int_0^t \int_{\partial\Omega} \alpha\frac{\partial G_\alpha}{\partial \vn_\vy}(\vx-\vy,t-\tau) u(\vy,\tau) \,\dif\vs_\vy \,\dif\tau, \\
  K_T' &\colon X' \to X', &K_T' w(\vx,t) &= \int_0^t \int_{\partial\Omega} \alpha\frac{\partial G_\alpha}{\partial \vn_\vx}(\vx-\vy,t-\tau) w(\vy,\tau) \,\dif\vs_\vy \,\dif\tau,
\end{align*}
where the integral representations on the right hold for sufficiently regular functions.

The above boundary integral equations are equivalent to the variational formulations
\begin{align}
  \label{eq:v_form_f_bie}
  \langle Vw, q \rangle_\Sigma &= \bigg\langle \bigg( \frac{1}{2}I + K \bigg) u, q \bigg\rangle_\Sigma &&\text{for all } q \in X',\\
  \label{eq:v_form_s_bie}
  \langle Du, r \rangle_\Sigma &= \bigg\langle \bigg( \frac{1}{2}I - K_T' \bigg) w, r \bigg\rangle_\Sigma &&\text{for all } r \in X
\end{align}
with the duality pairing $\langle \cdot, \cdot \rangle_\Sigma$ between $X'$ and $X$ given by the continuous extension of
\begin{equation*}
  \langle v, w \rangle_\Sigma := \int_0^T \int_{\partial\Omega} v(\vx,t) w(\vx,t) \,\dif\vs_\vx \,\dif t.
\end{equation*}
For the duality pairing with the hypersingular operator we have an alternative representation removing the non-integrable singularity, namely \cite{Cos1990, OfWat2021unpublished}
\begin{equation}
\label{eq:D_bil}
\langle Du, r \rangle_\Sigma = 
  \langle V \curl_{\partial\Omega} u ,\curl_{\partial\Omega} r \rangle_\Sigma
  - \langle \partial_t V (u \vn), r \vn \rangle_\Sigma,
\end{equation}
where the surface curl of a sufficiently regular function $u$ is defined by
\begin{equation*}
  \curl_{\partial\Omega} u(\vx,t) := \vn(x) \times \nabla_x \widetilde{u} (\vx, t),
\end{equation*}
for a suitable extension $\widetilde{u}$ of $u$ to an open neighbourhood of $\Sigma$. If the functions $u$ and $r$ are regular enough, \eqref{eq:D_bil} admits an integral representation. We will consider such a representation in the discrete setting in Section~\ref{subsec:03_hs}.

To solve an initial boundary value problem for the heat equation \eqref{eq:heat_eq} with homogeneous initial conditions and prescribed Dirichlet or Neumann boundary data, it suffices to determine the unknown boundary data. Then we can use the representation formula \eqref{eq:repr} to recover the solution. In the case of a Dirichlet boundary value problem the Neumann datum $w$ can be determined from~\eqref{eq:v_form_f_bie} while in the case of a Neumann boundary value problem the Dirichlet datum $u$ satisfies~\eqref{eq:v_form_s_bie}. It is shown in \cite{Cos1990, Noo1988} that these variational formulations admit a unique solution. In the next section we deal with their discretisation.

%% file: 03_bem.tex
For the discretisation of the variational formulations \eqref{eq:v_form_f_bie} and \eqref{eq:v_form_s_bie} we need a discretisation of the space-time boundary $\Sigma$. We restrict our attention to tensor product space-time discretisations~$\Sigma_h$ with uniform time steps. For a given uniform decomposition of the time interval
\begin{equation*}
  \overline{(0,T)} = \bigcup_{i=1}^{E_t} \overline{(t_{i-1}, t_i)} = \bigcup_{i=1}^{E_t} \overline{((i-1)h_t, i h_t)}
\end{equation*}
and an admissible triangular mesh $\Gamma_h$, which approximates $\Gamma := \partial\Omega$ and is given by
\begin{equation*}
  \Gamma_h = \bigcup_{j=1}^{E_\vx} \overline{\gamma_j}
\end{equation*}
with $\gamma_j$ denoting planar triangular elements, we define the space-time mesh
\begin{equation*}
  \overline{\Sigma_h} := \bigcup_{k=1}^{E_t E_{\vx}} \overline{\sigma_k} = \bigcup_{i=1}^{E_t} \bigcup_{j=1}^{E_\vx} \overline{\gamma_j} \times \overline{(t_{i-1}, t_i)}.
\end{equation*}
On $\Sigma_h$ we construct approximating spaces $X_h^{1,0} \subset X$ and $X_h^{0,0} \subset X'$ accordingly as tensor products, i.e. as linear combinations of functions whose spatial and temporal contributions can be separated as
\begin{equation*}
  \varphi_{\vx t,k}(\vx,t) = \varphi_{t,i} (t) \varphi_{\vx,j} (\vx).
\end{equation*}
We thus define the space
\begin{equation*}
  X_h^{0,0} := \vspan (\varphi^{0,0}_{\vx t,k})_{k=1}^{E_t E_{\vx}} = \vspan ((\varphi^0_{t,i}\varphi^0_{\vx,j})_{j=1}^{E_{\vx}})_{i=1}^{E_t}
\end{equation*}
of functions piecewise constant both in space and time and the space
\begin{equation*}
  X_h^{1,0} := \vspan (\varphi^{1,0}_{\vx t,k})_{k=1}^{E_t N_\vx} = \vspan ((\varphi^0_{t,i}\varphi^1_{\vx,j})_{j=1}^{N_{\vx}})_{i=1}^{E_t}
\end{equation*}
of functions globally continuous and piecewise linear in space and piecewise constant in time. Here we denote by $N_\vx$ the number of nodes of the triangular mesh $\Gamma_h$. 

A function $u_h$ in $X_h^{1,0}$ admits the representation 
\begin{equation*}
  u_h = \sum_{i=1}^{E_t} \sum_{j=1}^{N_\vx} u_{i,j} \varphi^0_{t,i} \varphi^1_{\vx,j},
\end{equation*}
where the first index of the coefficient $u_{i,j}$ is associated with time and the second one with space. This notation is slightly inconsistent with respect to the naming convention classically used for the function spaces, where the first superscript is related to space and the second one to time. However, for the implementation and representation of the matrices it is more natural to sort with respect to time first, which is why we use this notation.

To discretise the variational formulations \eqref{eq:v_form_f_bie} and \eqref{eq:v_form_s_bie} we replace the functions in $X$ and $X'$ with their discrete counterparts in $X_h^{1,0}$ and $X_h^{0,0}$ respectively. In the following subsections we give more details about the resulting discrete operators. In particular, we focus on the computation of the corresponding integrals.

\subsection{Single-layer boundary integral operator} \label{subsec:slo}

We start with the discretisation of the bilinear form $\langle V w, q \rangle_\Sigma$. By replacing $w$ with the approximation
\begin{equation*}
  w_h := \sum_{i=1}^{E_t} \sum_{j=1}^{E_\vx} w_{i,j} \varphi^0_{t,i} \varphi^0_{\vx,j}
\end{equation*}
and testing with a basis function
\begin{equation*}
  q_h := \varphi^0_{t,k} \varphi^0_{\vx,\ell}
\end{equation*}
we obtain
\begin{align*}
  \langle &Vw_h, q_h \rangle_{\Sigma_h} \\
  &= \int_0^T \int_{\Gamma_h} \varphi^0_{t,k}(t) \varphi^0_{\vx,\ell}(\vx) \int_0^t \int_{\Gamma_h} G_\alpha(\vx-\vy,t-\tau) \sum_{i=1}^{E_t} \sum_{j=1}^{E_\vx} w_{i,j} \varphi^0_{t,i}(\tau) \varphi^0_{\vx,j}(\vy) \,\dif\vs_\vy \,\dif\tau \,\dif\vs_\vx \,\dif t \\
  &= \sum_{i=1}^{k-1} \sum_{j=1}^{E_\vx} w_{i,j} \int_{\gamma_\ell} \int_{\gamma_j} \int_{t_{k-1}}^{t_k} \int_{t_{i-1}}^{t_i} G_\alpha(\vx-\vy,t-\tau) \,\dif\tau \,\dif t \,\dif\vs_\vy \,\dif\vs_\vx \\
  &+ \sum_{j=1}^{E_\vx} w_{k,j} \int_{\gamma_\ell} \int_{\gamma_j} \int_{t_{k-1}}^{t_k} \int_{t_{k-1}}^{t} G_\alpha(\vx-\vy,t-\tau) \,\dif\tau \,\dif t \,\dif\vs_\vy \,\dif\vs_\vx.
\end{align*}
Here we changed the order of the integrals, which is justified by Fubini's theorem and the fact that $G_\alpha$ is Lebesgue integrable on $\Sigma_h \times \Sigma_h$, which follows from the estimate \cite[Ch.~13 §3]{Pog1966}
\begin{equation} \label{eq:est_heat_ker}
  |G_\alpha(\vx-\vy,t-\tau)|\leq c(\alpha) \frac{1}{(t-\tau)^{3/4}} \frac{1}{\|\vx-\vy\|^{3/2}}.
\end{equation}

Since the fundamental solution depends only on the difference $t-\tau$ and the considered decomposition of the time interval  is uniform, the double temporal integrals depend only on the difference $d := k - i$. The duality pairing with all basis functions thus leads to the block Toeplitz matrix vector product
\begin{equation}
\label{eq:V_mat}
  \matV_h \vw =
  \begin{bmatrix}
    \matV_h^0 & 0 & \ldots & 0 \\
    \matV_h^1 & \ddots & \ddots & \vdots \\
    \vdots & \ddots & \ddots & 0 \\
    \matV_h^{E_t-1} & \dots & \matV_h^1 & \matV_h^0
  \end{bmatrix}
  \begin{bmatrix}
    \vw^0 \\[1mm]
    \vw^1 \\[1mm]
    \vdots \\[1mm]
    \vw^{E_t-1}
  \end{bmatrix}
\end{equation}
with vector components $w^d_j := w_{d+1,j}$ and the spatial matrix blocks defined by
\begin{equation}
\label{eq:V_blocks}
\begin{aligned}
  \matV_h^0[\ell,j] &:= \int_{\gamma_\ell} \int_{\gamma_j} \int_{0}^{t_1} \int_{0}^{t} G_\alpha(\vx-\vy,t-\tau) \,\dif\tau \,\dif t \,\dif\vs_\vy \,\dif\vs_\vx, \\
  \matV_h^d[\ell,j] &:= \int_{\gamma_\ell} \int_{\gamma_j} \int_{t_d}^{t_{d+1}} \int_{0}^{t_1} G_\alpha(\vx-\vy,t-\tau) \,\dif\tau \,\dif t \,\dif\vs_\vy \,\dif\vs_\vx
\end{aligned}
\end{equation}
for $d \in \{1,\ldots,E_t-1\}$. To set up $\matV_h$ we use analytic integration in time and a regularised quadrature in space as used in stationary problems. The details are given in the following paragraphs.
\subsubsection*{Temporal antiderivatives:}
Using $t_d = dh_t$, we have to evaluate
\begin{equation}
\label{eq:V_d}
  V^d(\vr) :=
  \begin{cases}
  \displaystyle
    \int_{0}^{h_t} \int_{0}^{t} G_\alpha(\vr,t-\tau) \,\dif\tau \,\dif t & \text{for } d = 0, \\[4mm]
  \displaystyle
    \int_{dh_t}^{(d+1)h_t} \int_{0}^{h_t} G_\alpha(\vr,t-\tau) \,\dif\tau \,\dif t & \text{for } d \in \{1,\ldots,E_t-1\}.
  \end{cases}
\end{equation}
We start with the latter. Integrating with respect to $\tau$ leads to
\begin{equation*}
  V^d(\vr) = \int_{dh_t}^{(d+1)h_t} \left( G^{\dif\tau}_\alpha(\vr,t-h_t) - G^{\dif\tau}_\alpha(\vr,t) \right) \,\dif t
\end{equation*}
with
\begin{equation}
  \label{eq:G_dtau}
  G^{\dif\tau}_\alpha(\vr,\delta) = \frac{1}{4\pi \alpha \|\vr\|} \erf\bigg( \frac{\|\vr\|}{2\sqrt{\alpha\delta}} \bigg)
\end{equation}
and the error function 
\begin{equation*}
  \erf(x) := \frac{2}{\sqrt{\pi}} \int_0^x e^{-t^2} \,\dif t.
\end{equation*}
Continuing in the integration we obtain
\begin{equation}
  \label{eq:vd}
  V^d(\vr) = 2G^{\dif\tau\dif t}_\alpha(\vr,d h_t) - G^{\dif\tau\dif t}_\alpha(\vr,(d+1)h_t) - G^{\dif\tau\dif t}_\alpha(\vr,(d-1)h_t)
\end{equation}
with 
\begin{equation}\label{eq:G_dtau_dt}
  G^{\dif\tau\dif t}_\alpha(\vr,\delta) = \frac{1}{4\pi} \bigg[ \bigg( \frac{\|\vr\|}{2\alpha^2} + \frac{\delta}{\alpha \|\vr\|} \bigg) \erf\bigg( \frac{\|\vr\|}{2\sqrt{\alpha\delta}} \bigg) + \frac{\sqrt{\delta}}{\sqrt{\pi\alpha^3}}\exp\bigg( -\frac{\|\vr\|^2}{4\alpha\delta} \bigg) \bigg]
\end{equation}
for $\delta > 0$ and $\|\vr\| > 0$.
For $V^0(\vr)$ we obtain in a similar fashion
\begin{align}\label{eq:V0}
  V^0(\vr) &= \int_{0}^{h_t} \left(G^{\dif \tau}_\alpha(\vr,0) - G^{\dif \tau}_\alpha(\vr, t)\right) \,\dif t 
  =h_t G^{\dif\tau}_\alpha(\vr,0) - G^{\dif\tau\dif t}_\alpha(\vr,h_t) + G^{\dif\tau\dif t}_\alpha(\vr,0).
\end{align}
Thus the integrals in \eqref{eq:V_blocks} are linear combinations of the integrals
\begin{align}
  \label{eq:V_spat_int_reg}
  &\int_{\gamma_\ell} \int_{\gamma_j} G^{\dif\tau\dif t}_\alpha(\vx - \vy, \delta) \,\dif\vs_\vy \,\dif\vs_\vx, \quad \text{for} \ \delta \in \{0, h_t, \dots, E_t h_t \}, \\
  \label{eq:V_spat_int_sing}
  &\int_{\gamma_\ell} \int_{\gamma_j} G^{\dif\tau}_\alpha(\vx - \vy, 0) \,\dif\vs_\vy \,\dif\vs_\vx.
\end{align}
Notice that a contribution with a fixed $\delta$ can be reused to assemble $\matV_h^d$ for several values of $d$.
\subsubsection*{Stable Evaluations of $V^d(\vr)$ for special cases:}
 We have to provide stable alternatives of \eqref{eq:G_dtau} and \eqref{eq:G_dtau_dt} for cases where the standard form does not allow an evaluation by a computer. In~\eqref{eq:vd} with $d=1$  and~\eqref{eq:V0} where $d=0$ we evaluate $G^{\dif\tau\dif t}_\alpha(\vr,\delta)$ from~\eqref{eq:G_dtau_dt} in $\delta = 0$ by
\begin{equation*}
  \lim_{\delta \to 0+} G^{\dif\tau\dif t}_\alpha(\vr,\delta) = \frac{\|\vr\|}{8\pi\alpha^2}
\end{equation*}
for $\|\vr\| > 0$. Similarly, we have to consider the limit
\begin{equation*}
  \lim_{\|\vr\|\to 0+} G^{\dif\tau\dif t}_\alpha(\vr,\delta) = \frac{\sqrt{\delta}}{2\sqrt{\pi^3\alpha^3}}
\end{equation*}
to evaluate  $G^{\dif\tau\dif t}_\alpha(\vr,\delta)$ in $\vr = 0$ for $\delta > 0$.
In~\eqref{eq:V0} we have to additionally evaluate  $G^{\dif\tau}_\alpha(\vr,\delta)$ from~\eqref{eq:G_dtau} in $\delta = 0$ by 
\begin{equation}
  \label{eq:G_dtau_lim_d}
  \lim_{\delta \to 0+} G^{\dif\tau}_\alpha(\vr,\delta) = \frac{1}{4\pi\alpha\|\vr\|}
\end{equation}
for $\|\vr\| > 0$.

\subsubsection*{Computation of the Galerkin weights of $\matV_h^d$:}
We have to compute the spatial integrals of \eqref{eq:V_spat_int_reg} and \eqref{eq:V_spat_int_sing}.
For $\delta > 0$ the integrand in \eqref{eq:V_spat_int_reg} is smooth. Therefore, standard quadrature routines can be applied to evaluate
\begin{equation}
\label{eq:quad_reg}
  \int_{\gamma_\ell} \int_{\gamma_j} G^{\dif\tau\dif t}_\alpha(\vx - \vy, \delta) \,\dif\vs_\vy \,\dif\vs_\vx 
  = 4 \Delta_\ell \Delta_j \int_{\hat\gamma} \int_{\hat\gamma} \widehat{G}^{\dif\tau\dif t}_\alpha(\hat\vx-\hat\vy, \delta) \,\dif\vs_{\hat\vy} \,\dif\vs_{\hat\vx}
\end{equation}
where we make use of the standard mapping to a reference element $\hat\gamma$. Here we denote the composition of the  mapping and the kernel function by  $\widehat{G}^{\dif\tau\dif t}_\alpha$ and the surface area of a triangle $\gamma_j$ by $\Delta_j$.

The integrand $G^{\dif\tau}_\alpha(\vx-\vy, 0)$ of \eqref{eq:V_spat_int_sing} given in~\eqref{eq:G_dtau_lim_d} has a singularity like the Laplace kernel for $\vx = \vy$. Thus we can use standard quadrature routines only if the triangles $\gamma_\ell$ and $\gamma_j$ are separated. If instead $\overline{\gamma_\ell}$ and~$\overline{\gamma_j}$ have non-empty intersection, i.e.~they share a vertex or an edge or are identical, we use regularised quadrature techniques based on the Duffy substitution~\cite{SauSch2010,ManTau2019}. The integral \eqref{eq:V_spat_int_sing} then transforms to an integral of the type
\begin{equation}
    \label{eq:quad_duffy}
  \sum_{s=1}^{N_{\mathrm{S}}} \int_0^1 \int_0^1 \int_0^1 \int_0^1 \widehat{G}^{\dif\tau}_\alpha(\vF^s_\vx(\eta_1,\eta_2,\eta_3,\xi)-\vF^s_\vy(\eta_1,\eta_2,\eta_3,\xi)) J^s(\eta_1,\eta_2,\eta_3,\xi) \,\dif\eta_1 \,\dif\eta_2 \,\dif\eta_3 \,\dif\xi
\end{equation}
with a mapping $\vF^s = (\vF^s_\vx, \vF^s_\vy)  \colon [0,1]^4 \to S \subset \hat\gamma \times \hat\gamma$ and the Jacobian $J^s \colon [0,1]^4 \to \R$,
\begin{equation*}
  \vF^s(\eta_1,\eta_2,\eta_3,\xi) = (\hat\vx, \hat\vy), \quad J^s(\eta_1,\eta_2,\eta_3,\xi) \,\dif\eta_1 \,\dif\eta_2 \,\dif\eta_3 \,\dif\xi = \dif\vs_\vy \,\dif\vs_\vx.
\end{equation*}

Analogously we deal with the integrals in \eqref{eq:V_spat_int_reg} for $\delta = 0$. Although in that case the function $\vr \mapsto G^{\dif\tau\dif t}_\alpha(\vr, 0)$ does not have a pole at $\vr=\vzero$ it is still not smooth and we use the regularised quadrature for intersecting triangles as well. This also unifies the implementation for other kernels possibly singular in this case.

If we use discrete test and trial functions with higher polynomial degree in space for the discretisation of the bilinear form $\langle V u, q \rangle$, e.g.~$u_h, q_h \in X_h^{1,0}$, the computation of the matrix entries follows the same lines. In particular, the matrix entries of the $d$-th block are given by
\begin{equation*}
  \matV_h^d[\ell,j] := \int_{\Gamma_h} \int_{\Gamma_h} \varphi_{\vx,\ell}(x) \varphi_{\vy,j}(y) V^d(\vx - \vy) \,\dif\vs_\vy \,\dif\vs_\vx,
\end{equation*}
with $V^d$ from \eqref{eq:V_d}.

\subsection{Double-layer boundary integral operator}\label{subsec:03_dl}
For the discretisation of $\langle K u, q \rangle_\Sigma$ we replace $u$ with its approximation $u_h$ in~$X_h^{1,0}$, i.e.
\begin{equation*}
  u_h := \sum_{i=1}^{E_t} \sum_{j=1}^{N_\vx} u_{i,j} \varphi^0_{t,i} \varphi^1_{\vx,j}.
\end{equation*}
By testing with the basis function
\begin{equation*}
  q_h := \varphi^0_{t,k} \varphi^0_{\vx,\ell}
\end{equation*}
we obtain for $d=k-i$ that
\begin{align*}
  \langle Ku_h, q_h \rangle_{\Sigma_h} &= \sum_{d=1}^{k-1} \sum_{j=1}^{N_\vx} u_{k-d,j} \int_{\gamma_\ell} \int_{\Gamma_h} \varphi^1_{\vx,j}(\vy) \alpha \int_{t_d}^{t_{d+1}} \int_{0}^{t_1} \frac{\partial G_\alpha}{\partial \vn_\vy}(\vx-\vy,t-\tau) \,\dif\tau \,\dif t \,\dif\vs_\vy \,\dif\vs_\vx \\
  &+ \sum_{j=1}^{N_\vx} u_{k,j} \int_{\gamma_\ell} \int_{\Gamma_h} \varphi^1_{\vx,j}(\vy) \alpha \int_{0}^{t_1} \int_{0}^{t} \frac{\partial G_\alpha}{\partial \vn_\vy}(\vx-\vy,t-\tau) \,\dif\tau \,\dif t \,\dif\vs_\vy \,\dif\vs_\vx.
\end{align*}
Again, we changed the order of integration using Fubini's theorem, which is applicable since $\partial G_\alpha/\partial \vn_\vy$ is integrable on~$\Sigma_h \times \Sigma_h$. This follows from the estimate~\cite[Ch.~13 §3]{Pog1966}
\begin{equation} \label{eq:est_n_deriv}
  \left|\frac{\partial G_\alpha}{\partial \vn_\vy}(\vx-\vy,t-\tau)\right| =
  \frac{|(\vx-\vy) \cdot \vn_\vy|}{16\pi^{3/2} (\alpha (t - \tau))^{5/2}} 
  \exp\bigg( -\frac{\|\vx-\vy\|^2}{4\alpha(t-\tau)} \bigg) 
  \leq c(\alpha) \frac{1}{(t-\tau)^{3/4}} \frac{|(\vx-\vy) \cdot \vn_\vy|}{\|\vx-\vy\|^{7/2}}
\end{equation}
similarly as in \cite[Sect.~8.2.2]{Hac1995}, because the spatial boundary $\Gamma_h$ is piecewise smooth. An analogous estimate to \eqref{eq:est_n_deriv} holds for $|\nabla_{\vy} G_\alpha|$. This allows us to exchange the spatial gradient and the time integrals for $\vx-\vy\neq 0$ \cite[Prop.~5.9]{BroKer2015} which yields
\begin{align*}
  \int_{0}^{t_1} \int_{0}^{t} \frac{\partial G_\alpha}{\partial \vn_\vy}(\vx-\vy,t-\tau) \,\dif\tau \,\dif t
  &= \frac{\partial}{\partial \vn_\vy} \int_{0}^{t_1} \int_{0}^{t} G_\alpha(\vx-\vy,t-\tau) \,\dif\tau \,\dif t, \\
  \int_{t_d}^{t_{d+1}} \int_{0}^{t_1} \frac{\partial G_\alpha}{\partial \vn_\vy}(\vx-\vy,t-\tau) \,\dif\tau \,\dif t
  &= \frac{\partial}{\partial \vn_\vy} \int_{t_d}^{t_{d+1}} \int_{0}^{t_1} G_\alpha(\vx-\vy,t-\tau) \,\dif\tau \,\dif t.
\end{align*}
\subsubsection*{Temporal antiderivatives:}
As before we analytically evaluate the integrals
\begin{equation*}
  K^d(\vr) :=
  \begin{cases}
  \displaystyle
    \alpha \frac{\partial}{\partial \vn_\vy} \int_{0}^{h_t} \int_{0}^{t} G_\alpha(\vr,t-\tau) \,\dif\tau \,\dif t & \text{for } d = 0, \\[4mm]
  \displaystyle
    \alpha \frac{\partial}{\partial \vn_\vy} \int_{dh_t}^{(d+1)h_t} \int_{0}^{h_t} G_\alpha(\vr,t-\tau) \,\dif\tau \,\dif t & \text{for } d \in \{1,\ldots,E_t-1\}.
  \end{cases}
\end{equation*}
For $d > 0$ we obtain from \eqref{eq:vd} that
\begin{align}\label{eq:K_d}
  K^d(\vr) = \alpha \bigg[ 2 \frac{\partial G^{\dif\tau\dif t}_\alpha}{\partial \vn_\vy} (\vr,dh_t) - \frac{\partial G^{\dif\tau\dif t}_\alpha}{\partial \vn_\vy}(\vr,(d+1)h_t) - \frac{\partial G^{\dif\tau\dif t}_\alpha}{\partial \vn_\vy}(\vr,(d-1)h_t) \bigg].
\end{align}
Since $G^{\dif\tau\dif t}_\alpha$ depends only on the norm of its first argument, we can write
\begin{equation}
  \label{eq:G_to_g}
  G^{\dif\tau\dif t}_\alpha(\vr,\delta) =: g_\alpha^{\dif\tau\dif t}(\|\vr\|,\delta) \quad\text{with}\quad g^{\dif\tau\dif t}_\alpha(\rho,\delta)\colon \R\times\R\to\R
\end{equation}
to get
\begin{equation}
\label{eq:dn_G_dt_dtau}
\begin{aligned}
  \frac{\partial G^{\dif\tau\dif t}_\alpha}{\partial \vn_\vy}(\vr,\delta) &= \frac{\partial g^{\dif\tau\dif t}_\alpha}{\partial \rho}(\|\vr\|,\delta) \vn_\vy \cdot \nabla_\vy \|\vx-\vy\|
  = - \frac{\partial g^{\dif\tau\dif t}_\alpha}{\partial \rho}(\|\vr\|,\delta) \frac{ \vr \cdot \vn_\vy}{\| \vr\|}
\end{aligned}
\end{equation}
with
\begin{equation*}
  \frac{\partial g^{\dif\tau\dif t}_\alpha}{\partial \rho}(\|\vr\|,\delta) =\frac{1}{4\pi} \bigg[  \bigg( \frac{1}{2\alpha^2} - \frac{\delta}{\alpha\|\vr\|^2} \bigg) \erf\bigg( \frac{\|\vr\|}{2\sqrt{\alpha\delta}} \bigg) + \frac{\sqrt{\delta}}{\|\vr\|\sqrt{\pi\alpha^3}} \exp \bigg( -\frac{\|\vr\|^2}{4\alpha\delta} \bigg) \bigg].
\end{equation*}
Collecting all intermediate steps brings us to
\begin{align}\label{eq:dG}
  \alpha\frac{\partial G^{\dif\tau\dif t}_\alpha}{\partial \vn_\vy}(\vr,\delta) &= -\frac{1}{4\pi} \frac{ \vr \cdot \vn_\vy}{\|\vr\|}\bigg[  \bigg( \frac{1}{2\alpha} - \frac{\delta}{\|\vr\|^2} \bigg) \erf\bigg( \frac{\|\vr\|}{2\sqrt{\alpha\delta}} \bigg)
  + \frac{\sqrt{\delta}}{\|\vr\|\sqrt{\pi\alpha}} \exp \bigg( -\frac{\|\vr\|^2}{4\alpha\delta} \bigg) \bigg].
\end{align}

For $K^0(\vr)$ we use \eqref{eq:V0} to get
\begin{align}
  K^0(\vr) = \alpha \bigg[ h_t\frac{\partial G^{\dif\tau}_\alpha}{\partial \vn_\vy} (\vr,0) - \frac{\partial G^{\dif\tau\dif t}_\alpha}{\partial \vn_\vy}(\vr,h_t) + \frac{\partial G^{\dif\tau\dif t}_\alpha}{\partial \vn_\vy}(\vr,0) \bigg]. \label{eq:K0}
\end{align}
Similarly as in \eqref{eq:dn_G_dt_dtau} with \eqref{eq:G_to_g} we have
\begin{equation*}
  \frac{\partial G^{\dif\tau}_\alpha}{\partial \vn_\vy}(\vr,\delta) = -\frac{\partial g^{\dif\tau}_\alpha}{\partial \rho}(\|\vr\|,\delta) \frac{\vr \cdot \vn_\vy}{\|\vr\|}
\end{equation*}
with $G^{\dif\tau}_\alpha(\vr,\delta) =: g^{\dif\tau}_\alpha(\|\vr\|,\delta)$ and
\begin{equation*}
  \frac{\partial g^{\dif\tau}_\alpha}{\partial \rho}(\|\vr\|,\delta) = -\frac{1}{4\pi} \frac{1}{\|\vr\|}\bigg[ \frac{1}{\alpha\|\vr\|} \erf\bigg( \frac{\|\vr\|}{2\sqrt{\alpha\delta}} \bigg) - \frac{1}{\sqrt{\pi\alpha^3\delta}} \exp\bigg( -\frac{\|\vr\|^2}{4\alpha\delta} \bigg) \bigg].
\end{equation*}
Thus, we obtain
\begin{equation}
\label{eq:dn_G_dtau}
  \alpha \frac{\partial G^{\dif\tau}_\alpha}{\partial \vn_\vy}(\vr,\delta) = \frac{1}{4\pi} \frac{\vr \cdot \vn_\vy}{\|\vr\|^2}\bigg[ \frac{1}{\|\vr\|} \erf\bigg( \frac{\|\vr\|}{2\sqrt{\alpha\delta}} \bigg) - \frac{1}{\sqrt{\pi\alpha\delta}} \exp\bigg( -\frac{\|\vr\|^2}{4\alpha\delta} \bigg) \bigg].
\end{equation}

\subsubsection*{Stable Evaluations of $K^d(\vr)$ for special cases:}

For a stable evaluation of \eqref{eq:K_d} for $d=1$, we evaluate~\eqref{eq:dG} for $\delta = 0$ and $\|\vr\|>0$ by
\begin{equation*}
  \lim_{\delta \to 0_+}  \alpha\frac{\partial G^{\dif\tau\dif t}_\alpha}{\partial \vn_\vy}(\vr,\delta)  = -\frac{\vr \cdot \vn_\vy}{8 \pi\alpha \| \vr\|}.
\end{equation*}
Conversely, for $\delta > 0$ the value of $\alpha(\partial G^{\dif\tau\dif t}_\alpha/\partial \vn_\vy)(\vr,\delta)$ in $\vr = 0$ is given by
\begin{equation*}
  \lim_{\|\vr\| \to 0_+} \alpha\frac{\partial G^{\dif\tau\dif t}_\alpha}{\partial \vn_\vy}(\vr,\delta) = 0.
\end{equation*}
This follows by estimating
\begin{align*}
	\left| \alpha\frac{\partial G_\alpha^{\dif\tau \dif t}}{\partial \vn_\vy}
	(\vr, \delta) \right|
	&\leq \left|
	\frac{\left(\frac{\| \vr \|^2}{2\alpha} - \delta \right)
	\erf\left(\frac{\| \vr \|}{2\sqrt{\alpha\delta}}  \right)
	+ \frac{\sqrt{\delta} \| \vr \|}{\sqrt{\pi\alpha}}
	\exp \left( -\frac{\| \vr\|^2}{4\alpha \delta} \right)}
	{4 \pi \| \vr\|^2 }
	 \right|
	 =: \big|\tilde g^{\dif \tau \dif t}_\alpha(\|\vr \|, \delta)\big|
\end{align*}
and observing (e.g. using L'Hospital's rule) that the limit of $\tilde g^{\dif \tau \dif t}_\alpha$ for $\rho = \| \vr\| \rightarrow 0_+$ is zero. 

For a stable evaluation of $K^0(\vr)$ by~\eqref{eq:K0}, we provide the values of~\eqref{eq:dn_G_dtau}
in $\delta = 0$ for $\|\vr\| > 0$ as
\begin{equation}
\label{eq:dn_G_dtau_lim_r}
  \lim_{\delta \to 0_+} \alpha\frac{\partial G^{\dif\tau}_\alpha}{\partial \vn_\vy}(\vr,\delta) = \frac{ \vr \cdot \vn_\vy}{4\pi\|\vr\|^3}.
\end{equation}

\subsubsection*{Computation of the Galerkin weights of $\matK_h^d$:}
The layout of the block matrix $\matK_h$ is the same as the one of $\matV_h$ in \eqref{eq:V_mat}, i.e. 
\begin{equation}
\label{eq:K_mat}
  \matK_h =
  \begin{bmatrix}
    \matK_h^0 & 0 & \ldots & 0 \\
    \matK_h^1 & \ddots & \ddots & \vdots \\
    \vdots & \ddots & \ddots & 0 \\
    \matK_h^{E_t-1} & \dots & \matK_h^1 & \matK_h^0
  \end{bmatrix}.
\end{equation}
The individual blocks $\matK_h^d$ are set up as
\begin{equation*}
  \matK^d_h[ \ell , j ] = \int_{\gamma_\ell} \int_{\Gamma_h} \varphi^1_{\vx,j}(\vy) K^d(\vx-\vy) \,\dif\vs_\vy \,\dif\vs_\vx,
\end{equation*}
where the integrals are evaluated in the same way as the integrals we considered for the single-layer operator. In particular, we use the same regularisation technique. This time, we have to deal with a singularity similar to the one of the double-layer boundary integral operator of the Laplacian, see \eqref{eq:dn_G_dtau_lim_r}.

\begin{remark*}[The Galerkin matrix of the operator $K_T'$]
The matrix $\matK_h^{\top_\vx}$ related to the operator~$K_T'$ is obtained from $\matK_h$ by transposing each block in~\eqref{eq:K_mat} separately, not the matrix as a whole.
\end{remark*}

\subsection{Hypersingular boundary integral operator} \label{subsec:03_hs}
For functions $u_h$ and $r_h$ in $X_h^{1,0}$ one can show that the right-hand side of \eqref{eq:D_bil} and therefore the bilinear form $\langle D u_h, r_h \rangle_{\Sigma_h}$ admits the following weakly singular integral representation~(see~\cite{OfWat2021unpublished}, compare also \cite{Mes2014,DohNiiSte2019})
\begin{equation} \label{eq:D_bil_int_rep}
\begin{aligned}
  &\langle Du_h, r_h \rangle_{\Sigma_h} \\
  &= \alpha^2 \int_0^T \int_{\Gamma_h} \curl_{\Gamma_h} r_h(\vx,t) \cdot \bigg( \int_0^t \int_{\Gamma_h} \curl_{\Gamma_h} u_h(\vy,\tau) G_\alpha(\vx-\vy,t-\tau) \,\dif\vs_\vy \,\dif\tau \bigg) \,\dif\vs_\vx \,\dif t \\
  &- \alpha \sum_{n=1}^{E_t} \int_{t_{n-1}}^{t_n} \int_{\Gamma_h} r_h(\vx,t) \vn(\vx) \cdot \bigg[ \int_0^{t_{n-1}} \int_{\Gamma_h} \vn(\vy) u_h(\vy,\tau) \frac{\partial G_\alpha}{\partial \tau}(\vx-\vy,t-\tau) \,\dif\vs_\vy \,\dif\tau \\
  &\phantom{- \alpha \sum_{n=1}^{E_t} \int_{t_{n-1}}^{t_n} \int_{\Gamma_h} r_h(\vx,t) \vn(\vx) \cdot \bigg[ } - \int_{\Gamma_h} \vn(\vy) u_h(\vy,t_{n-1}{+}) G_\alpha(\vx-\vy,t-t_{n-1}) \, \dif\vs_\vy \bigg] \,\dif\vs_\vx \,\dif t,
\end{aligned}
\end{equation}
where $u_h(\vy,t_{n-1}{+})$ denotes the right limit of $u_h$ with respect to time in $t_{n-1}$. 

By inserting
\begin{equation*}
  u_h := \sum_{i=1}^{E_t} \sum_{j=1}^{N_\vx} u_{i,j} \varphi^0_{t,i} \varphi^1_{\vx,j}
\end{equation*}
and testing with the basis function
\begin{equation*}
  r_h := \varphi^0_{t,k} \varphi^1_{\vx,\ell}
\end{equation*}
we obtain for $d=k-i$ that
\begin{equation} \label{eq:D_bil_int_rep_split}
\begin{aligned}
  &\langle Du_h, r_h \rangle_{\Sigma_h} = \\ 
  &\sum_{d=1}^{k-1} \sum_{j=1}^{N_\vx} u_{k-d,j}\int_{\Gamma_h} \int_{\Gamma_h} \!\! \curl_{\Gamma_h} \varphi^1_{\vx,\ell}(\vx) \cdot \curl_{\Gamma_h} \varphi^1_{\vx,j}(\vy) \alpha^2 \int_{t_{d}}^{t_{d+1}} \!\!\! \int_{0}^{t_1} \! G_\alpha(\vx-\vy,t-\tau) \,\dif\tau \,\dif t \,\dif\vs_\vy \,\dif\vs_\vx \\
    &+\sum_{j=1}^{N_\vx} u_{k,j}\int_{\Gamma_h} \int_{\Gamma_h} \curl_{\Gamma_h} \varphi^1_{\vx,\ell}(\vx) \cdot \curl_{\Gamma_h} \varphi^1_{\vx,j}(\vy) \alpha^2 \int_{0}^{t_1} \int_{0}^{t} G_\alpha(\vx-\vy,t-\tau) \,\dif\tau \,\dif t \,\dif\vs_\vy \,\dif\vs_\vx \\
  &- \sum_{d=1}^{k-1} \sum_{j=1}^{N_\vx} u_{k-d,j} \int_{\Gamma_h} \int_{\Gamma_h} \vn(\vx) \cdot \vn(\vy) \varphi^1_{\vx,\ell}(\vx) \varphi^1_{\vx,j}(\vy) \alpha \int_{t_d}^{t_{d+1}} \!\!\! \int_0^{t_1} \frac{\partial G_\alpha}{\partial \tau}(\vx-\vy,t-\tau) \,\dif\tau \,\dif t \,\dif\vs_\vy \,\dif\vs_\vx \\
    &+ \sum_{j=1}^{N_\vx} u_{k,j} \int_{\Gamma_h} \int_{\Gamma_h} \vn(\vx) \cdot \vn(\vy) \varphi^1_{\vx,\ell}(\vx) \varphi^1_{\vx,j}(\vy) \alpha \int_{0}^{t_1} G_\alpha(\vx-\vy,t) \,\dif t \,\dif\vs_\vy \,\dif\vs_\vx.
\end{aligned}
\end{equation}
Changing the order of the integrals is justified as before by Fubini's theorem. Indeed, for the integrals in the first two lines we can argue as in the case of the single-layer operator in Section~\ref{subsec:slo}. For the integrals in the fourth line it suffices to observe that $G_\alpha$ is Lebesgue integrable on $\Sigma_h \times \Gamma_h$ which follows from \eqref{eq:est_heat_ker}. Similarly, for the integrals in the third line we note that $\partial G_\alpha/ \partial\tau$ is Lebesgue integrable on all sets $(\Gamma_h \times (t_d, t_{d+1}))\times(\Gamma_h \times (0, t_{1}))$ with $d \geq 1$, because there holds the estimate
\begin{align*}
  \left|\frac{\partial G_\alpha}{\partial \tau}(\vx-\vy,t-\tau)\right| =
  \frac{|6\alpha(t-\tau) - \|\vx-\vy\|^2|}{(4 \alpha)^{5/2} \pi^{3/2} (t - \tau)^{7/2}} 
  \exp\bigg( -\frac{\|\vx-\vy\|^2}{4\alpha(t-\tau)} \bigg) 
  \leq c(\alpha) \frac{1}{(t-\tau)^{7/4}} \frac{1}{\|\vx-\vy\|^{3/2}}.
\end{align*}
\subsubsection*{Temporal antiderivatives:}

For the first summands in \eqref{eq:D_bil_int_rep_split} we know the temporal antiderivatives from the single-layer boundary integral operator, see~\eqref{eq:vd} and~\eqref{eq:V0}. 

For the second part in \eqref{eq:D_bil_int_rep_split} we analytically evaluate the integrals
\begin{equation*}
  D^{2,d}(\vr) :=
  \begin{cases}
  \displaystyle
    \alpha \int_{0}^{h_t} G_\alpha(\vr,t) \,\dif t & \text{for } d = 0, \\[4mm]
  \displaystyle
    -\alpha \int_{dh_t}^{(d+1)h_t} \int_{0}^{h_t} \frac{\partial G_\alpha}{\partial \tau}(\vr,t-\tau) \,\dif\tau \,\dif t & \text{for } d \in \{1,\ldots,E_t-1\}.
  \end{cases}
\end{equation*}
For $d>0$ we can write
\begin{align}
  D^{2,d}(\vr) &= -\alpha \int_{dh_t}^{(d+1)h_t} G_\alpha(\vr,t-h_t) - G_\alpha(\vr,t) \,\dif t \nonumber \\
  &= -\alpha \big[ 2 G_\alpha^{\dif t}( \vr, dh_t ) - G_\alpha^{\dif t}(\vr, (d + 1)h_t) - G_\alpha^{\dif t}(\vr, (d - 1)h_t ) \big]\label{eq:D2d}
\end{align}
with $G_\alpha^{\dif t} = -G_\alpha^{\dif \tau}$ from \eqref{eq:G_dtau}.
For $d=0$ we directly get
\begin{equation*}
  D^{2,0}(\vr) = -\alpha \big[ G_\alpha^{\dif t}( \vr, 0 ) - G_\alpha^{\dif t}(\vr, h_t ) \big].
\end{equation*}
\subsubsection*{Stable Evaluations of $D^{2,d}(\vr)$ for special cases:}
The values $G_\alpha^{\dif t}(\vr, 0)$ in \eqref{eq:D2d} for $\vr \neq 0$ are obtained by the limit in~\eqref{eq:G_dtau_lim_d}. Additionally we have to treat the limit for $\delta>0$,
\begin{equation*}
  \lim_{\|\vr\| \to 0_+} G^{\dif t}_\alpha(\vr,\delta) = -\frac{1}{4\sqrt{\pi^3\alpha^3\delta}}
\end{equation*}
to evaluate $G_\alpha^{\dif t}(\vr, \delta)$ in $\vr = 0$ in a stable way.

\subsubsection*{Galerkin matrix $\matD_h$:}
The Galerkin matrix $\matD_h$ possesses the same layout as \eqref{eq:V_mat} and can be split into $\matD_h = \matD_h^1 + \matD_h^2$. For $D_h^2$ we have just computed the temporal antiderivatives. Its blocks~$\matD_h^{2,d}$ are set up as
\begin{equation*}
  \matD^{2,d}_h[ \ell , j ] = \int_{\Gamma_h} \int_{\Gamma_h} \vn(\vx) \cdot \vn(\vy) \varphi^1_{\vx,\ell}(\vx) \varphi^1_{\vx,j}(\vy) D^{2,d}(\vx-\vy) \,\dif\vs_\vy \,\dif\vs_\vx.
\end{equation*}
Again, these integrals are handled in the same way as those of the single-layer operator.

For the matrix $\matD_h^1$ emerging from the first two lines of \eqref{eq:D_bil_int_rep_split} we can make use of $\matV_h$. Since the surface curls of spatially piecewise linear functions are piecewise constant on triangles, we can rewrite any of the summands in the first part of \eqref{eq:D_bil_int_rep_split} as
\begin{align*}
  \matD^{1,d}_h[\ell,j] &:= \int_{\Gamma_h} \int_{\Gamma_h}\curl_{\Gamma_h} \varphi^1_{\vx,\ell}(\vx) \cdot \curl_{\Gamma_h} \varphi^1_{\vx,j}(\vy) \alpha^2 \int_{t_{d}}^{t_{d+1}} \int_{0}^{t_1} G_\alpha(\vx-\vy,t-\tau) \,\dif\tau \,\dif t \,\dif\vs_\vy \,\dif\vs_\vx \\
  &= \sum_{\gamma_n \subset \supp \varphi^1_{\vx,\ell}} \sum_{\gamma_m \subset \supp \varphi^1_{\vx,j}} \curl_{\Gamma_h} \varphi^1_{\vx,\ell}|_{\gamma_n}(\vx) \cdot \curl_{\Gamma_h} \varphi^1_{\vx,j}|_{\gamma_m}(\vy) \alpha^2 \matV_h^d[n,m].
\end{align*}
Thus, for all $d \in \{0,\ldots,E_t-1\}$ the block $\matD_h^{1,d}$ is a sparse transformation of the single-layer block~$\matV^d_h$ from \eqref{eq:V_blocks}. In particular, the individual blocks can be assembled by
\begin{equation} \label{eq:D_1_sp_trafo}
  \matD_h^{1,d} = \matT^\top
  \begin{bmatrix}
    \alpha^2\matV^d_h & \matO & \matO \\
    \matO & \alpha^2\matV^d_h & \matO \\
    \matO & \matO & \alpha^2\matV^d_h
  \end{bmatrix}
  \matT,\quad
  \matT :=
  \begin{bmatrix}
    \matT_1 \\
    \matT_2 \\
    \matT_3
  \end{bmatrix},\quad
  \matT_o[m,j] := [\curl_{\Gamma_h} \varphi^1_{\vx,j}|_{\gamma_m}]_o.
\end{equation}

\subsection{Boundary integral equations and systems of linear equations}

To solve Dirichlet initial boundary value problems, we consider a Galerkin variational formulation of the weakly singular boundary integral equation~\eqref{eq:v_form_f_bie} with $u=g$ and end up with the system of linear equations
\begin{equation*}
  \matV_h \vw = \bigg(\frac{1}{2} \matM_h + \matK_h \bigg) \vg
\end{equation*}
with an $L_2(\Sigma_h)$ projection of the Dirichlet data $g$ into $X_h^{1,0}$. 
We have described the matrices $\matV_h$ and $\matK_h$ in Sections~\ref{subsec:slo}-\ref{subsec:03_dl}. 
The mass matrix $\matM_h$ realises the identity in \eqref{eq:v_form_f_bie} and has the form
\begin{equation*}
  \matM_h = h_t
    \begin{bmatrix}
      \matM_{h,\vx} & 0 & \ldots & 0 \\
      0 & \matM_{h,\vx} & \ddots & \vdots \\
      \vdots & \ddots & \ddots & 0 \\
      0 & \dots &0 & \matM_{h,\vx}
    \end{bmatrix}, \quad
  \matM_{h,\vx}[\ell, j] 
    = \int_{\Gamma_h} \varphi^0_{\vx,\ell}(\vx) \varphi^1_{\vx,j}(\vx)\,\dif\vs_\vx.
\end{equation*}

For a Neumann initial boundary value problem we solve the Galerkin variational formulation of the hypersingular boundary integral equation~\eqref{eq:v_form_s_bie} with $w=h$. The related system of linear equations is
\begin{equation*}
  \matD_h \vu = \bigg(\frac{1}{2} \matM_h^{\top_\vx} - \matK_h^{\top_\vx} \bigg) \vh
\end{equation*}  
with an $L_2(\Sigma_h)$ projection of the Neumann data $h$ into $X_h^{0,0}$.
We have presented the matrix~$\matD_h$ in Section~\ref{subsec:03_hs}. In addition we need to assemble the matrices $\matK_h$ and $\matM_h$. The blockwise transposition can be realised in the application of the matrices.

\subsection{Single- and double-layer potential} \label{subsec:slp_dlp}

To evaluate the discretised representation formula \eqref{eq:repr} in $\vx\in\Omega$ and $t_k + \varepsilon = k h_t + \varepsilon$ with $\varepsilon \in [0, h_t)$ we have to compute the contribution of the single-layer potential
\begingroup
\allowdisplaybreaks
\begin{align*}
 \widetilde{V}&w_h(\vx, t_k + \varepsilon) = \sum_{i=1}^{k} \sum_{j=1}^{E_\vx} w_{i,j} \int_{\gamma_j} \int_{t_{i-1}}^{t_i} G_\alpha(\vx-\vy,t_k + \varepsilon -\tau) \,\dif\tau \,\dif\vs_\vy \\
 &+ \sum_{j=1}^{E_\vx} w_{k+1,j} \int_{\gamma_j} \int_{t_k}^{t_k + \varepsilon} G_\alpha(\vx-\vy,t_k + \varepsilon -\tau) \,\dif\tau \,\dif\vs_\vy \\
 &= \sum_{d=0}^{k-1} \sum_{j=1}^{E_\vx} w_{k-d,j} \int_{\gamma_j} \int_{0}^{h_t} G_\alpha(\vx-\vy,(d+1)h_t + \varepsilon-\tau) \,\dif\tau \,\dif\vs_\vy \\
 &+ \sum_{j=1}^{E_\vx} w_{k+1,j} \int_{\gamma_j} \int_{0}^{\varepsilon} G_\alpha(\vx-\vy,\varepsilon-\tau) \,\dif\tau \,\dif\vs_\vy \\
  &= \sum_{j=1}^{E_\vx} \int_{\gamma_j} \sum_{d=0}^{k-1} w_{k-d,j} \big[ G_\alpha^{\dif\tau}(\vx-\vy, dh_t + \varepsilon) - G_\alpha^{\dif\tau}(\vx-\vy, (d+1)h_t + \varepsilon) \big] \,\dif\vs_\vy \\
  &+ \sum_{j=1}^{E_\vx} \int_{\gamma_j} w_{k+1,j} \big[ G_\alpha^{\dif\tau}(\vx-\vy, 0) - G_\alpha^{\dif\tau}(\vx-\vy, \varepsilon ) \big] \,\dif\vs_\vy
\end{align*}
\endgroup
with the antiderivative $G_\alpha^{\dif\tau}$ known from \eqref{eq:G_dtau} and its limit for $\delta \to 0_+$ in \eqref{eq:G_dtau_lim_d}. Since we evaluate the potential in points $\vx$ which are not on $\Gamma_h$, all integrands are smooth and standard quadrature can be used to compute all integrals.

For the double-layer potential we similarly obtain
\begin{align*}
 W&u_h(\vx, t_k + \varepsilon) = \sum_{i=1}^{k} \sum_{j=1}^{N_\vx} u_{i,j} \int_{\Gamma_h} \varphi^1_{\vx,j}(\vy) \int_{t_{i-1}}^{t_i} \alpha \frac{\partial G_\alpha}{\partial \vn_\vy}(\vx-\vy,t_k + \varepsilon-\tau) \,\dif\tau \,\dif\vs_\vy \\
 &+ \sum_{j=1}^{N_\vx} u_{k+1,j} \int_{\Gamma_h} \varphi^1_{\vx,j}(\vy) \int_{t_k}^{t_k + \varepsilon} \alpha \frac{\partial G_\alpha}{\partial \vn_\vy}(\vx-\vy,t_k + \varepsilon-\tau) 
 \,\dif\tau \,\dif\vs_\vy \\
  &= \sum_{j=1}^{N_\vx}\int_{\Gamma_h} \varphi^1_{\vx,j}(\vy) \sum_{d=0}^{k-1} u_{k-d,j} \alpha\bigg[ \frac{\partial G_\alpha^{\dif\tau}}{\partial \vn_\vy}(\vx-\vy, d h_t + \varepsilon) - \frac{\partial G_\alpha^{\dif\tau}}{\partial \vn_\vy}(\vx-\vy, (d+1)h_t + \varepsilon) \bigg] \,\dif\vs_\vy \\
  &+ \sum_{j=1}^{N_\vx} \int_{\Gamma_h} \varphi^1_{\vx,j}(\vy) u_{k+1,j} \alpha\bigg[ \frac{\partial G_\alpha^{\dif\tau}}{\partial \vn_\vy}(\vx-\vy, 0) - \frac{\partial G_\alpha^{\dif\tau}}{\partial \vn_\vy}(\vx-\vy, \varepsilon) \bigg] \,\dif\vs_\vy
\end{align*}
with $\alpha\partial G_\alpha^{\dif\tau} / \partial \vn_\vy$ from \eqref{eq:dn_G_dtau} and its limit for $\delta \to 0_+$ in \eqref{eq:dn_G_dtau_lim_r}.

%% file: 035_tldr.tex
For a better readability and to provide a reference to a reader implementing the method we provide a summary of the developed formulae below.

With uniform time steps all matrices $\matA_h \in \{\matV_h, \matK_h, \matK_h^{\top_\vx}, \matD_h\}$ possess a block Toeplitz structure

\begin{equation*}
  \matA_h =
  \begin{bmatrix}
    \matA_h^0 & 0 & \ldots & 0 \\
    \matA_h^1 & \ddots & \ddots & \vdots \\
    \vdots & \ddots & \ddots & 0 \\
    \matA_h^{E_t-1} & \dots & \matA_h^1 & \matA_h^0
  \end{bmatrix}.
\end{equation*}
The hypersingular operator matrix is built as $\matD_h^d = \matD_h^{1,d} + \matD_h^{2,d}$ with
\begin{equation*}
  \matD_h^{1,d} = \matT^\top
  \begin{bmatrix}
    \alpha^2\matV^d_h & \matO & \matO \\
    \matO & \alpha^2\matV^d_h & \matO \\
    \matO & \matO & \alpha^2\matV^d_h
  \end{bmatrix}
  \matT,\quad
  \matT :=
  \begin{bmatrix}
    \matT_1 \\
    \matT_2 \\
    \matT_3
  \end{bmatrix},\quad
  \matT_o[m,j] := [\curl_{\partial\Omega} \varphi^1_{\vx,j}|_{\gamma_m}]_o.
\end{equation*}
The matrix $\matK_h^{\top_\vx}$ discretising the operator $K_T'$ is obtained from $\matK_h$ by blockwise transposition.
Individual blocks $\matA_h^d \in \{\matV_h^d, \matK_h^d, \matD_h^{2,d}\}$ are built by a standard regularised BEM quadrature as
\begin{equation*}
  \matA_h^d[\ell,j] = \int_{\Gamma_h} \int_{\Gamma_h} A^d(\vx-\vy) \varphi_\ell(\vx) \varphi_j(\vy) \,\dif\vs_\vy \,\dif\vs_\vx
\end{equation*}
with $A^d \in \{V^d, K^d, D^{2,d}\}$ and
\begin{align*}
  V^d(\vr) &= 2G^{\dif\tau\dif t}_\alpha(\vr,d h_t) - G^{\dif\tau\dif t}_\alpha(\vr,(d+1)h_t) - G^{\dif\tau\dif t}_\alpha(\vr,(d-1)h_t), \\
  V^0(\vr) &= h_t G^{\dif\tau}_\alpha(\vr,0) - G^{\dif\tau\dif t}_\alpha(\vr,h_t) + G^{\dif\tau\dif t}_\alpha(\vr,0), \\[4mm]
  K^d(\vr) &= \alpha \bigg[ 2 \frac{\partial G^{\dif\tau\dif t}_\alpha}{\partial \vn_\vy} (\vr,dh_t) - \frac{\partial G^{\dif\tau\dif t}_\alpha}{\partial \vn_\vy}(\vr,(d+1)h_t) - \frac{\partial G^{\dif\tau\dif t}_\alpha}{\partial \vn_\vy}(\vr,(d-1)h_t) \bigg], \\
  K^0(\vr) &= \alpha \bigg[ h_t\frac{\partial G^{\dif\tau}_\alpha}{\partial \vn_\vy} (\vr,0) - \frac{\partial G^{\dif\tau\dif t}_\alpha}{\partial \vn_\vy}(\vr,h_t) + \frac{\partial G^{\dif\tau\dif t}_\alpha}{\partial \vn_\vy}(\vr,0) \bigg], \\[4mm]
  D^{2,d}(\vr) &= -\alpha \big[ 2 G_\alpha^{\dif t}( \vr, dh_t ) - G_\alpha^{\dif t}(\vr, (d + 1)h_t) - G_\alpha^{\dif t}(\vr, (d - 1)h_t ) \big], \\
  D^{2,0}(\vr) &= -\alpha \big[ G_\alpha^{\dif t}( \vr, 0 ) - G_\alpha^{\dif t}(\vr, h_t ) \big].
\end{align*}
The antiderivatives of the heat kernel and limit cases for stable evaluations are given by
\begingroup
\allowdisplaybreaks
\begin{align*}
  G^{\dif\tau\dif t}_\alpha(\vr,\delta) &= \frac{1}{4\pi} \bigg[ \bigg( \frac{\|\vr\|}{2\alpha^2} + \frac{\delta}{\alpha \|\vr\|} \bigg) \erf\bigg( \frac{\|\vr\|}{2\sqrt{\alpha\delta}} \bigg) + \frac{\sqrt{\delta}}{\sqrt{\pi\alpha^3}}\exp\bigg( -\frac{\|\vr\|^2}{4\alpha\delta} \bigg) \bigg], \\
  \lim_{\delta \to 0+} G^{\dif\tau\dif t}_\alpha(\vr,\delta) &= \frac{\|\vr\|}{8\pi\alpha^2} \quad \text{for } \|\vr\| > 0, \\
  \lim_{\|\vr\|\to 0+} G^{\dif\tau\dif t}_\alpha(\vr,\delta) &= \frac{\sqrt{\delta}}{2\sqrt{\pi^3\alpha^3}} \quad \text{for } \delta > 0, \\[4mm]
  G^{\dif\tau}_\alpha(\vr,\delta) &= \frac{1}{4\pi \alpha \|\vr\|} \erf\bigg( \frac{\|\vr\|}{2\sqrt{\alpha\delta}} \bigg), \\
  \lim_{\delta \to 0+} G^{\dif\tau}_\alpha(\vr,\delta) &= \frac{1}{4\pi\alpha\|\vr\|} \quad \text{for } \|\vr\| > 0, \\
  \alpha\frac{\partial G^{\dif\tau\dif t}_\alpha}{\partial \vn_\vy}(\vr,\delta) &= -\frac{1}{4\pi} \frac{ \vr \cdot \vn_\vy}{\|\vr\|}\bigg[  \bigg( \frac{1}{2\alpha} - \frac{\delta}{\|\vr\|^2} \bigg) \erf\bigg( \frac{\|\vr\|}{2\sqrt{\alpha\delta}} \bigg)
  + \frac{\sqrt{\delta}}{\|\vr\|\sqrt{\pi\alpha}} \exp \bigg( -\frac{\|\vr\|^2}{4\alpha\delta} \bigg) \bigg], \\
  \lim_{\delta \to 0_+}  \alpha\frac{\partial G^{\dif\tau\dif t}_\alpha}{\partial \vn_\vy}(\vr,\delta) &= -\frac{\vr \cdot \vn_\vy}{8 \pi\alpha \| \vr\|} \quad \text{for } \|\vr\| > 0, \\
  \lim_{\|\vr\| \to 0_+} \alpha\frac{\partial G^{\dif\tau\dif t}_\alpha}{\partial \vn_\vy}(\vr,\delta) &= 0 \quad \text{for } \delta > 0, \\[4mm]
  \alpha \frac{\partial G^{\dif\tau}_\alpha}{\partial \vn_\vy}(\vr,\delta) &= \frac{1}{4\pi} \frac{\vr \cdot \vn_\vy}{\|\vr\|^2}\bigg[ \frac{1}{\|\vr\|} \erf\bigg( \frac{\|\vr\|}{2\sqrt{\alpha\delta}} \bigg) - \frac{1}{\sqrt{\pi\alpha\delta}} \exp\bigg( -\frac{\|\vr\|^2}{4\alpha\delta} \bigg) \bigg], \\
  \lim_{\delta \to 0_+} \alpha\frac{\partial G^{\dif\tau}_\alpha}{\partial \vn_\vy}(\vr,\delta) &= \frac{ \vr \cdot \vn_\vy}{4\pi\|\vr\|^3} \quad \text{for } \|\vr\| > 0, \\[4mm]
G_\alpha^{\dif t}(\vr,\delta) &= -\frac{1}{4\pi \alpha \|\vr\|} \erf\bigg( \frac{\|\vr\|}{2\sqrt{\alpha\delta}} \bigg), \\
  \lim_{\|\vr\| \to 0_+} G^{\dif t}_\alpha(\vr,\delta) &= -\frac{1}{4\sqrt{\pi^3\alpha^3\delta}} \quad \text{for } \delta > 0.
\end{align*}
\endgroup

%% file: 04_imp.tex
In this section we discuss an implementation strategy for the assembly of the single-layer matrix $\matV_h$. All other BEM matrices can be treated analogously. The computationally most intensive part is the evaluation of the antiderivatives $G^{\dif\tau}$ and $G^{\dif\tau\dif t}$. Indeed, in addition to the evaluation of the distance between spatial coordinates $\|\vx-\vy\|$, which is the most time consuming part of the BEM assembly for the Laplace equation, one has to evaluate the exponential and error functions in many quadrature points for all blocks of the Toeplitz matrix. The implementation strategy in shared memory thus follows the ideas presented by the authors previously for 3d space and 2d space-time BEM in~\cite{ZapMerMal2016,ZapOfMer2018,DohZapOfMer2019}. To make use of modern multicore processors with vector arithmetic units we make use of features of modern OpenMP~\cite{omp5}, namely threading and SIMD vectorisation. The source code of the library besthea implemented by the authors is publicly available~\cite{besthea}.

\subsection{Assembly of blocks}

The naive approach to assemble the Toeplitz matrix \eqref{eq:V_mat} would be to assemble blocks $\matV^d_h$ one by one, i.e. loop over the parameter~$d$. Looking at~\eqref{eq:vd}, this would mean that $G^{\dif\tau\dif t}(\cdot,\delta)$ would have to be evaluated multiple times in all spatial quadrature points for a fixed $\delta$ and different values of~$d$. E.g., for $\delta = 2h_t$ the same kernel would have to be evaluated for all blocks with~$d \in \{1,2,3\}$. 

Taking into account the uniform discretisation of the time interval one can instead loop over $i_t := \delta/h_t$ and thus evaluate the costly kernel once only. In Tables~\ref{tab:deltad} and \ref{tab:ddelta} we summarise the relation between $d$ and $i_t$, i.e.~we state to which blocks $\matV^d_h$ the kernels $G^{\dif\tau\dif t}(\cdot,i_t h_t)$ contribute and vice versa. Note that a similar strategy can also be applied to evaluate the single- and double-layer potentials given in Section~\ref{subsec:slp_dlp}.

\begin{table}[ht]
\centering
\begin{tabular}{r|r|r|r|r|r|r}
$i_t = \delta/h_t$ & $0$ & $1$ & $2$ & $\cdots$ & $E_t - 1$ & $E_t$\\
\midrule
$d$ & $0$, $1$  & $0$, $1$, $2$ & $1$, $2$, $3$ & $\cdots$ & $E_t-2$, $E_t-1$ & $E_t-1$
\end{tabular}
\caption{Mapping of variables $i_t \to d$}
\label{tab:deltad}
\bigskip
\begin{tabular}{r|r|r|r|r|r}
$d$ & $0$ & $1$ & $2$ & $\cdots$ & $E_t-1$\\
\midrule
$i_t = \delta/h_t$ & $0$, $1$  & $0$, $1$, $2$ & $1$, $2$, $3$ & $\cdots$ & $E_t-2$, $E_t-1$, $E_t$
\end{tabular}
\caption{Mapping of variables $d \to i_t$}
\label{tab:ddelta}
\end{table}

A sketch of the matrix assembly code is given in Listing~\ref{list:assembly}. As pointed out above, we loop over the variable $i_t = \delta/h_t$ and continue with visiting all test and trial spatial elements (triangles). The loop over test triangles is distributed among available OpenMP threads in a dynamic fashion. For each pair of elements the functions {\lstinline[basicstyle=\ttfamily]+evaluate_kernel+} and {\lstinline[basicstyle=\ttfamily]+add_to_matrix+} are called to assemble the local contribution and add it to the global matrix multiplied with the test and trial basis functions. We give more details about these procedures in the next subsection.

\begin{lstlisting}[label={list:assembly}, caption={General structure of the matrix assembling.}]
for (int i_t = 0; i_t <= n_timesteps; ++i_t) {
  ...
|*#pragma omp for schedule(dynamic)*|
  for (int i_test = 0; i_test < n_elements; ++i_test) {
    ...
    for (int i_trial = 0; i_trial < n_elements; ++i_trial) {
      ...
      evaluate_kernel(i_test, i_trial, i_t * h_t, ...);
      add_to_matrix(i_test, i_trial, i_t, ...);
} } }
\end{lstlisting}

\subsection{Local contributions}

To exploit the full potential of floating point units we vectorise the code at the level of local contributions to the global matrix. For simplicity we opt for the OpenMP implementation of vector processing similarly as in \cite{ZapMerMal2016,ZapOfMer2018,DohZapOfMer2019}.

Looking back on~\eqref{eq:quad_duffy}, we approximate the regularised integrals of the type
\begin{equation}
   \label{eq:ss}
    \int_0^1 \int_0^1 \int_0^1 \int_0^1 f(\eta_1, \eta_2, \eta_3, \xi) \,\dif\eta_1 \,\dif\eta_2 \,\dif\eta_3 \,\dif\xi \approx \sum_{i,j,k,\ell=1}^M w_i w_j w_k w_\ell f(z_i, z_j, z_k, z_\ell).
\end{equation}
by a tensor product quadrature scheme defined in $[0,1]^4$. The regular integrals can be evaluated by triangle rules \cite[Section C1]{RjaSte2007} as
\begin{equation}
   \label{eq:st}
    \int_\gamma \int_\gamma f(\vx,\vy) \,\dif\vy \,\dif\vx \approx \sum_{i,j=1}^N w_i w_j f(\vz_i, \vz_j).
\end{equation}
\begin{lstlisting}[label={list:kernel}, caption={{\lstinline[basicstyle=\ttfamily]+evaluate_kernel+}.}]
|*#pragma omp simd aligned(x1, x2, x3, y1, y2, y3, kernel, w : 64) simdlen(8)*|
for (int k = 0; k < size; ++k) {
  kernel[k] =
   _kernel->eval(x1[k] - y1[k], x2[k] - y2[k], x3[k] - y3[k], delta, ...) * w[k]; 
} 
\end{lstlisting}
In both cases we collapse the sums, or loop, into a single one to make the vector of quadrature points as long as possible to evaluate the kernel function efficiently. This is shown in Listing~\ref{list:kernel}, where the temporal antiderivative is evaluated. The variable {\lstinline[basicstyle=\ttfamily]+size+} corresponds to $M^4$ and $N^2$ from~\eqref{eq:ss} and \eqref{eq:st}, respectively. The OpenMP pragma tells the compiler that SIMD vectorisation should be used, that all the underlying arrays are aligned at the 64-byte boundary, and that the vector size should be 8 (we assume double precision arithmetic and AVX512 instruction set extension). Notice that we also make use of the structure of arrays concept separating coordinates of the quadrature nodes into separate arrays {\lstinline[basicstyle=\ttfamily]+x1+}, {\lstinline[basicstyle=\ttfamily]+x2+}, {\lstinline[basicstyle=\ttfamily]+x3+}, \dots, to ensure unit strided access to data. Earlier work \cite{ZapMerMal2016,ZapOfMer2018,DohZapOfMer2019} has shown that this approach is more efficient than an array of vectors.

After performing and storing the kernel evaluations in  {\lstinline[basicstyle=\ttfamily]+kernel+} by {\lstinline[basicstyle=\ttfamily]+evaluate_kernel+} for the current pair of elements, we evaluate the test and trial basis functions as shown in Listing~\ref{list:global}, multiply with {\lstinline[basicstyle=\ttfamily]+kernel+} and add {\lstinline[basicstyle=\ttfamily]+value+} to the respective spatial and temporal indices in the global matrix. Here we use the mapping from Table~\ref{tab:deltad}. The multiplier for {\lstinline[basicstyle=\ttfamily]+value+} is determined from~\eqref{eq:vd}. Again we make use of vectorisation. The {\lstinline[basicstyle=\ttfamily]+add_atomic+} function makes use of the OpenMP atomic clause to avoid data races between individual threads.

\begin{lstlisting}[label={list:global}, caption={{\lstinline[basicstyle=\ttfamily]+add_to_matrix+}.}]
for (int i = 0; i < n_loc_test; ++i) {
  for (int j = 0; j < n_loc_trial; ++j) {
    value = 0.0;
    
|*#pragma omp simd aligned(x1_ref, x2_ref, y1_ref, y2_ref, kernel : 64) \*|
|*private(test, trial) reduction(+ : value) simdlen(8)*|
    for (long k = 0; k < size; ++k) {
      test = test_basis.eval(x1_ref[k], x2_ref[k], ...);
      trial = trial_basis.eval(y1_ref[k], y2_ref[k], ...);
      value += kernel[k] * test * trial;
    }
    value *= test_area * trial_area;
      
    if (i_t > 0) {
      matrix.add_atomic(i_t - 1, test_l2g[i], trial_l2g[j], -value);
      if (i_t < n_timesteps) {
        matrix.add_atomic(i_t, test_l2g[i], trial_l2g[j], 2.0 * value);
      }
    } else {
      matrix.add_atomic(0, test_l2g[i], trial_l2g[j], value );
    }
    if (i_t < n_timesteps - 1) {
      matrix.add_atomic(i_t + 1, test_l2g[i], trial_l2g[j], -value);
} } }
\end{lstlisting}

%% file: 05_exp.tex
In this section we perform numerical experiments validating the presented approach both in terms of convergence and scalability in shared memory. The experiments have been performed at the Barbora supercomputer at IT4Innovations National Supercomputing Center, Czech Republic.

\subsection{Convergence}

First of all, we check that the presented semi-analytic evaluation of the integrals and its implementation in~\cite{besthea} is correct. To that end we consider the initial problem \eqref{eq:heat_eq}--\eqref{eq:initial:cond} with the heat capacity constant $\alpha = 0.5$ and zero initial conditions in the space-time domain ${Q := (-1,1)^3 \times (0,1)}$. We choose the solution $u(\vx,t) = G_{\alpha}(\vx - \vy^\ast, t)$ with $\vy^\ast := (0,0,1.5)^\top$, which allows us to validate our numerical approximation.
We consider both the Dirichlet problem with the prescribed boundary datum
\begin{equation*}
  u(\vx,t) = G_{\alpha}(\vx - \vy^\ast, t) \quad\text{for }(\vx,t) \in \Sigma,
\end{equation*}
and the Neumann problem with
\begin{equation*}
  \alpha \frac{\partial u}{\partial \vn} = \alpha \frac{\partial G_{\alpha}}{\partial \vn_\vx}(\vx - \vy^\ast, t) \quad\text{for }(\vx,t) \in \Sigma.
\end{equation*}
We have described details of the applied Galerkin methods in Section~\ref{sec:03_bem}. In addition we use the matrix $\matV_h^{11}$ as the essential part of an operator preconditioner for~$\matD_h$ \cite{SteWen1998,Hip2006,DohEtAl2018}, where $\matV_h^{11}$ is the realisation of the single-layer operator for functions piecewise linear and globally continuous in space.
We check the convergence of the approximations $u_h$ and $w_h$ corresponding to $\vu$ and $\vw$ to the known Cauchy data on a sequence of uniformly refined meshes $\Sigma_h$.

We consider only tensor product meshes $\Sigma_h$ in this paper. The coarsest one is formed by a surface mesh consisting of 192 triangular elements, i.e.~32 congruent triangles on each face of the cube, and a partition of the time interval $(0,1)$ into 8~time steps. At each refinement level we quadrisect all triangles and bisect the time steps, i.e.~we keep $h_\vx \approx h_t$. The solution of the BEM system is computed by the FGMRES~\cite{Saa1993} method with a relative accuracy of~$10^{-8}$. 

In Tables~\ref{tab:dirconv},~\ref{tab:neuconv} we provide the convergence results. In the first two columns, $E_t$ and~$E_\vx$ denote the number of elements in time and space, respectively. The columns labelled with $L^2(\Sigma_h)$ contain the relative errors
\begin{equation*}
  L^2(\Sigma_h)(u_h) := \frac{\|u-u_h\|_{L^2(\Sigma_h)}}{\|u\|_{L^2(\Sigma_h)}}
\end{equation*}
with
\begin{equation*}
  \|u\|^2_{L^2(\Sigma_h)} := \int_0^t \int_{\Gamma_h} |u(\vx, t)|^2 \,\dif\vs_\vy \,\dif\tau.
\end{equation*}
These integrals are evaluated using standard tensor product quadrature rules in space and time of sufficiently high orders. The estimated order of convergence provided in the columns denoted by eoc is computed as
\begin{equation*}
  \mathrm{eoc}(u_h) = \log_2\bigg( \frac{L^2(\Sigma_h)(u_{2h})}{L^2(\Sigma_h)(u_h)} \bigg).
\end{equation*}
For comparison, Tables~\ref{tab:dirconv} and~\ref{tab:neuconv} contain not only the results for the computed approximations $u_h$ and $w_h$, but also for the $L^2(\Sigma_h)$ projections $u_h^\ast$ and $w_h^\ast$ of the known solution defined by
\begin{equation*}
  u_h^\ast = \argmin_{v_h \in X_h^{1,0}} \|v_h - u\|_{L^2(\Sigma_h)}, 
  \quad w_h^\ast = \argmin_{z_h \in X_h^{0,0}} \|z_h - w\|_{L^2(\Sigma_h)}.
\end{equation*}
In the last two columns we present convergence results for the evaluation of the representation formula. For this purpose the representation formula was evaluated in $10^4$ nodes $(\widetilde \vx_j, \widetilde t_j)$ distributed in $[-0.5,0.5]^3 \times [0.25, 0.75]$. The columns labelled with~$\ell^2$ contain the relative errors
\begin{equation*}
  \ell^2(u_h) := \frac{\sqrt{\sum_{j} |u(\widetilde{\vx}_j, \widetilde{t}_j)-u_h(\widetilde{\vx}_j, \widetilde{t}_j)|^2}}
    {\sqrt{\sum_{j} |u(\widetilde{\vx}_j, \widetilde{t}_j)|^2}}.
\end{equation*}
  
\begin{table}[b]
\centering
\begin{tabular}{rr|rr|rr|rr}
 && \multicolumn{2}{c|}{computed $w_h$} & \multicolumn{2}{c|}{projected $w_h^\ast$} & \multicolumn{2}{c}{representation} \\
 \midrule
\multicolumn{1}{c}{$E_t$} & \multicolumn{1}{c|}{$E_\vx$} & \multicolumn{1}{c}{$L^2(\Sigma_h)$} & \multicolumn{1}{c|}{eoc} & \multicolumn{1}{c}{$L^2(\Sigma_h)$} & \multicolumn{1}{c|}{eoc} & \multicolumn{1}{c}{$\ell^2$} & \multicolumn{1}{c}{eoc} \\
\midrule
8 & 192 & 6.07e-1 & --- & 5.49e-1 & --- & 2.99e-2 & --- \\
16 & 768 & 4.28e-1 & 0.50 & 3.77e-1 & 0.54 & 3.46e-3 & 3.11 \\
32 & 3072 & 1.80e-1 & 1.25 & 1.70e-1 & 1.15 & 6.51e-4 & 2.41\\
64 & 12288 & 9.94e-2 & 0.86 & 9.58e-2 & 0.82 & 1.48e-4 & 2.14 \\
\end{tabular}
\caption{Dirichlet problem and the convergence of Neumann data.}
\label{tab:dirconv}
\end{table}

\begin{table}
\centering
\begin{tabular}{rr|rr|rr|rr}
 && \multicolumn{2}{c|}{computed $u_h$} & \multicolumn{2}{c|}{projected $u_h^\ast$} & \multicolumn{2}{c}{representation} \\
 \midrule
\multicolumn{1}{c}{$E_t$} & \multicolumn{1}{c|}{$E_\vx$} & \multicolumn{1}{c}{$L^2(\Sigma_h)$} & \multicolumn{1}{c|}{eoc} & \multicolumn{1}{c}{$L^2(\Sigma_h)$} & \multicolumn{1}{c|}{eoc} & \multicolumn{1}{c}{$\ell^2$} & \multicolumn{1}{c}{eoc} \\
\midrule
8 & 192 & 3.14e-1 & --- & 2.50e-1 & --- & 5.16e-2 & ---\\
16 & 768 & 1.51e-1 & 1.06 & 1.27e-1 & 0.98 & 1.57e-2 & 1.72 \\
32 & 3072 & 6.88e-1 & 1.13 & 6.11e-2 & 1.05 & 3.80e-3 & 2.04 \\
64 & 12288 & 3.45e-2 & 0.99 & 3.18e-2 & 0.94 & 1.09e-3 & 1.80 \\
\end{tabular}
\caption{Neumann problem and the convergence of Dirichlet data.}
\label{tab:neuconv}
\end{table}

Let us shortly comment on the results in Tables~\ref{tab:dirconv} and~\ref{tab:neuconv}. In both tables, the $L^2$ error of the computed approximation follows the best possible error, which is attained by the respective projections. In the case of the Dirichlet problem, the estimated orders of convergence of the $L^2$~errors vary quite a lot. Asymptotically we would expect at least an order of 0.75, while previous examples indicated that an order of 1 can be attained \cite[Thm.~7.4 and Sect.~8.2]{DohNiiSte2019}. Even though we are probably still in a preasymptotic regime due to the relatively small number of unknowns which we consider limited by the use of a standard, non-compressed BEM, our computations agree with these expectations. For the evaluation error inside the domain we expect and observe a quadratic convergence order \cite[Eq.~(7.5)]{DohNiiSte2019}. Also in the case of the Neumann problem our results agree with the theory. We expect and observe convergence order 1 for the $L^2$~error of the Dirichlet datum \cite[Eq.~(7.16)]{Doh2019} and order~1.5 for the evaluation error \cite[Sect.~7.2.2]{Doh2019}. A different refinement strategy of two subdivisioning steps in time with one spatial refinement step would provide better convergence rates.
In total, we observe expected convergence behaviours, which indicate the correctness of the developed and implemented quadrature routines.

\subsection{Scalability}
\begin{table}[b]
\centering
\begin{tabular}{r|rrrr|rrrr}
 & \multicolumn{4}{c|}{time [s]} & \multicolumn{4}{c}{efficiency [\%]} \\
\midrule
threads & \multicolumn{1}{c}{$\matV_h$} & \multicolumn{1}{c}{$\matV_h^{11}$} & \multicolumn{1}{c}{$\matK_h$} & \multicolumn{1}{c|}{$\matD_h$} & \multicolumn{1}{c}{$\matV_h$} & \multicolumn{1}{c}{$\matV_h^{11}$} & \multicolumn{1}{c}{$\matK_h$} & \multicolumn{1}{c}{$\matD_h$} \\
\midrule
1 & 119.62 & 308.03 & 118.89 & 443.45 & 100.00 & 100.00 & 100.00 & 100.00 \\
2 & 63.54 & 165.81 & 60.78 & 235.00 & 94.14 & 94.35 & 97.81 & 92.89 \\
4 & 31.09 & 88.97 & 31.74 & 124.32 & 96.20 & 89.18 & 93.65 & 86.55 \\
8 & 15.68 & 46.98 & 16.75 & 64.80 & 95.34 & 85.54 & 88.72 & 81.95 \\
16 & 8.06 & 24.64 & 8.26 & 33.48 & 92.80 & 82.78 & 90.00 & 78.13 \\ 
18 & 7.10 & 22.09 & 7.26 & 30.03 & 93.65 & 82.05 & 90.98 & 77.47 \\ 
\midrule
36 & 4.54 & 12.94 & 4.97 & 19.18 & 73.19 & 64.21 & 66.40 & 66.12 \\
\bottomrule
\end{tabular}
\caption{Scalability of the assembly of BEM matrices.}
\label{tab:asssc}
\end{table}

The scalability of the besthea library~\cite{besthea} has been tested on the same example as before, but on a fixed mesh with 3072 spatial boundary elements, 32 time steps, and the representation formula was evaluated in $1089 \cdot 32 = 34848$ space-time points. The library and examples were compiled by the Intel Compiler 19.0.5.281 with the flags {\lstinline[basicstyle=\ttfamily]+-O3 -qopenmp -xcore-avx512 -qopt-zmm-usage=high+} to make use of the AVX512 instruction set available on the 18-core Intel Xeon Gold 6240 CPU at the Barbora supercomputer. The nodes are configured as dual socket, i.e. every node consists of two such~CPUs.

The baseline for our experiments is given by the performance on a single thread. The number of threads is controlled by the {\lstinline[basicstyle=\ttfamily]+KMP_HW_SUBSET+} environment variable. When using up to 18 threads (a single socket) we set it to {\lstinline[basicstyle=\ttfamily]+KMP_HW_SUBSET=1s,Xc+} with {\lstinline[basicstyle=\ttfamily]+X+} denoting the number of threads. This ensures that the threads stay within a single socket. To use all 36 threads we set {\lstinline[basicstyle=\ttfamily]+KMP_HW_SUBSET=2s,18c+}. 

\begin{table}
\centering
\begin{tabular}{r|rr|rr}
 & \multicolumn{2}{c|}{time [s]} & \multicolumn{2}{c}{efficiency [\%]} \\
\midrule
threads & \multicolumn{1}{c}{$\widetilde{V}w_h$} & \multicolumn{1}{c|}{$Wu_h$} & \multicolumn{1}{c}{$\widetilde{V}w_h$} & \multicolumn{1}{c}{$W_hu_h$}\\
\midrule
1 & 8.36 & 21.17 & 100.00 & 100.00 \\
2 & 5.44 & 11.65 & 76.84 & 90.87 \\
4 & 2.40 & 5.68 & 87.08 & 93.14 \\
8 & 1.28 & 2.97 & 81.43 & 89.11 \\
16 & 0.69 & 1.53 & 76.09 & 86.49 \\
18 & 0.50 & 1.26 & 92.89 & 93.11 \\
\midrule
36 & 0.31 & 0.82 & 74.11 & 71.73 \\
\bottomrule
\end{tabular}
\caption{Scalability of the evaluation of potentials.}
\label{tab:reprsc}
\end{table}

In Table~\ref{tab:asssc} we provide the assembly times of BEM matrices and the efficiency of the code. 
One can see that the efficiency stays above 90 \% within a single socket with piecewise constant basis functions, where the contribution to the global matrix does not require atomic addition. For piecewise linear functions there are memory conflicts between individual threads and the efficiency is reduced, although the numbers stay reasonable. When crossing the socket and utilising all 36~cores the assembly times are further reduced, although the efficiency drops. This is caused by the fact that the matrix data is stored in {\lstinline[basicstyle=\ttfamily]+std::vector+} with the standard allocator which is not NUMA aware. The first touch policy thus cannot be easily applied and threads can access memory across sockets. This effect could be alleviated by a different storage structure or a raw array. However, the future aim of the besthea library is to use the fast multipole method parallelised in distributed memory via MPI and assign a single process per socket. 

Similarly, Table~\ref{tab:reprsc} provides scalability results for the evaluation of the single- and double-layer potentials, see Section~\ref{subsec:slp_dlp}. Again, the efficiency is above 90 \% when the whole socket is populated and drops when accessing both sockets.

%% file: 06_con.tex
The aim of the paper was to provide the readers with semi-analytical formulae for the assembly of boundary element matrices and the evaluation of the representation formula for the heat equation in three spatial dimensions. Throughout the paper a uniform discretisation of the timeline is chosen for simplicity, however, the same antiderivatives can be used on non-uniform grids. Moreover, the approach has been implemented in the publicly available \cpp{} library besthea~\cite{besthea} and thus our results can be used in further BEM projects. In the numerical experiments we have validated that the formulae deliver the expected results.

The provided implementation supplies a fast computation of the entries of the Galerkin matrices using threading and vectorisation. However it does not scale optimally across NUMA nodes (sockets). This is due to the fact that the main aim of the library is to provide boundary element methods accelerated by the fast multipole method (FMM) and parallelised in distributed memory. The formulae provided here will be used for the near field entries only, thus further optimisation of the full assembly is not planned (also taking into account the massive memory requirements). With FMM, a single MPI process will be assigned to a single socket and thus the NUMA effects will be automatically overcome.